\documentclass[12pt]{amsart}
\usepackage{latexsym, amsmath, amscd, amssymb, amsthm,longtable} 
\usepackage[english]{babel}
\newtheorem{te}{Theorem}
\newtheorem{pr}{Example}[section]

\usepackage{youngtab}
 
\usepackage{ytableau,tikz,varwidth}
\usepackage{scalerel}
\usepackage{cancel}
\usepackage{graphicx}

\usepackage{enumitem}

\begin{document}

\noindent

\title {  The  $\mathfrak{sl}_2$-actions on the symmetric polynomials and on  Young diagrams. } 

\author{Leonid Bedratyuk}
\address{Khmelnytskyi National University, Insituts'ka  st., 11,  Khmelnytskyi, 29016, Ukraine}
\begin{abstract}
In the article, two implementations of the representation of the complex Lie algebra $\mathfrak{sl}_2$ on the algebra of symmetric polynomials $\Lambda_n$ by differential operators are proposed. The realizations of irreducible subrepresentations, both finite-dimensional and infinite-dimensional, are described, and the decomposition of $\Lambda_n$ is found. The actions on  the  Schur polynomials is also determined. By using an  isomorphism between $\Lambda_n$ and the vector space of Young diagrams $\mathbb{Q}\mathcal{Y}_n$ with no more than $n$ rows, these representations are transferred to $\mathbb{Q}\mathcal{Y}_n$.

\end{abstract}
\maketitle

\section{Introduction}
Let $\Lambda_n$ be the algebra of symmetric polynomials in $n$ variables. We consider $\Lambda_n$ as an infinite-dimensional $\mathbb{Q}$-vector space on which we intend to introduce an $\mathfrak{sl}_2$-module structure, where $\mathfrak{sl}_2$ is the \textit{complex} Lie algebra with the basis $\langle \mathfrak{e}_-, \mathfrak{h}, \mathfrak{e}_+\rangle$ and the following commutation relations:
$$
[h,\mathfrak{e}_+]=2 \mathfrak{e}_+, [h,\mathfrak{e}_-]=-2 \mathfrak{e}_-, [\mathfrak{e}_+,\mathfrak{e}_-]=\mathfrak{h}.
$$
Recall that a representation $(\rho, V)$ of the algebra $\mathfrak{sl}_2$ in a vector space $V$ is a homomorphism $\rho: \mathfrak{sl}_2 \to {\rm End}(V)$ that preserves commutators. In this paper, we consider two ways of defining such a representation by differential operators. The first representation is defined by the following differential operators on $\Lambda_n$:
\begin{align*}
&\rho_1(\mathfrak{e}_+)=x_1^2\partial_{1}+x_2^2\partial_{2}+\cdots+x_n^2\partial_{n},\\
&\rho_1(\mathfrak{h})=2(x_1\partial_{1}+x_2\partial_{2}+\cdots+x_n\partial_{n}),\\
&\rho_1(\mathfrak{e}_-)=-(\partial_{1}+\partial_{2}+\cdots+\partial_{n}),
\end{align*}
where $\partial_i=\dfrac{\partial}{\partial x_i}$.

In the  paper we prove  that as an $\mathfrak{sl}_2$-module, $\Lambda_n$ is a direct sum of infinite-dimensional irreducible \textit{lowest weight modules} $W_\omega$ of $\mathfrak{sl}_2$, with a real weight $\omega > 0$, corresponding to type $([ \circ)$ in the notations of \cite{HT}. We explicitly describe all symmetric polynomials  which are lowest weight vectors of $\rho_1$ and prove that they, together with the polynomial $x_1+x_2+\cdots+x_n$, form a basis of the algebra $\Lambda_n$. To obtain this result, we extensively used results related to the theory of locally nilpotent derivations.

 For the case of two variables and natural $m$, as a by-product, we obtained the combinatorial identities:
\begin{align*}
&\sum_{k=0}^{2m-1} \left(-\frac{1}{2}\right)^{k+1} {2m+1 \choose k} (x^{2m+1-k}+y^{2m+1-k}) (x+y)^k  =\frac{m}{2^{2m}}(x+y)^{2m+1},\\
& \sum_{k=0}^{2m-2} (-1)^k 2^{2m-k-1} {2m\choose k} (x^{2m-k}+y^{2m-k}) (x+y)^k  =(2m-1) (x+y)^{2m}+(x-y)^{2m},
\end{align*}
which are similar to Girard-Waring power sum identities, see for example \cite{G}.

Also, we explicitly find the action of the operators of $\rho_1$ on the Schur polynomials $\boldsymbol{s}_\lambda$.

The second representation of $\mathfrak{sl}_2$ is defined by the following differential operators:
\begin{align*}
&\rho_2(\mathfrak{e}_+)=\sum_{k=1}^n(- x_k^2\partial_{k}+d \cdot x_k),\\ 
&\rho_2(\mathfrak{h})=2 \sum_{k=1}^n x_k\partial_{k}-n \cdot d, \\
&\rho_2(\mathfrak{e}_-) =\sum_{k=1}^n \partial_{k}.
\end{align*}
and depends on the parameter $d \in \mathbb{N}$. In $\Lambda_d$, there exists a finite-dimensional submodule $\Lambda_n^{(d)}$ which is isomorphic to the symmetric power ${\rm Sym}^n(V_d)$ where $V_d$ is the standard irreducible $\mathfrak{sl}_2$-module of dimension $d+1$. This isomorphism establishes a connection between symmetric polynomials and classical invariant theory.

In the last section, the obtained actions are transferred to the $\mathbb{Q}$-vector space $\mathbb{Q}\mathcal{Y}_n$ generated by  Young diagrams, which contain no more than $n$ rows. Interest in this question is inspired by the fact that Kerov operators, which define an $\mathfrak{sl}_2$-action on unbounded Young diagrams, no longer define an $\mathfrak{sl}_2$-action on $\mathcal{Y}_n$. Using the $\mathfrak{sl}_2$-isomorphism of vector spaces $\mathbb{Q}\mathcal{Y}_n$  and $\Lambda_n$, the results obtained for symmetric polynomials are transferred to  $\mathbb{Q}\mathcal{Y}_n$.

\vspace{1cm}

\section{The first $\mathfrak{sl}_2$-action on $\Lambda_n$. }

Consider the standard basis of the complex matrix Lie algebra $\mathfrak{sl}_2$:

$$
\mathfrak{e}_+=\begin{pmatrix} 0 & 1 \\ 0 & 0 \end{pmatrix}, \mathfrak{h}=\begin{pmatrix} 1 & 0 \\ 0 & -1 \end{pmatrix}, \mathfrak{e}_-=\begin{pmatrix} 0 & 0 \\ 1 & 0 \end{pmatrix},
$$

with the following commutation relations
$$
[h,\mathfrak{e}_+]=2 \mathfrak{e}+, [h,\mathfrak{e}_-]=-2 \mathfrak{e}_-, [\mathfrak{e}_+,\mathfrak{e}_-]=\mathfrak{h}.
$$
A representation $(\rho, V)$ of the algebra $\mathfrak{sl}_2$ in a vector space $V$ is called a homomorphism $\rho: \mathfrak{sl}_2 \to {\rm End}(V)$ that preserves the commutation relations. In the module terminology, the vector space $V$ is called an $\mathfrak{sl}_2$-module, and the homomorphism $\rho$ is called an $\mathfrak{sl}_2$-action.

The classification of $\mathfrak{sl}_2$-modules, both finite-dimensional and infinite dimensional, is given in \cite{HT}.
We will be interested in the irreducible infinite-dimensional \textit{lowest weight } $\mathfrak{sl}_2$-module $W_\omega,$ with real \textit{weight} $\omega > 0$, (of type $([ \circ)$ in the notation of \cite[Proposition 1.2.6]{HT}), which has a basis $v_0,v_1, \ldots, v_k, \ldots $ with the following action
\begin{align*}
&\rho(\mathfrak{e}_-) v_0=0,\\
&\rho(\mathfrak{e}_-) v_k=-k(\omega+k-1) v_{k-1},\\
&\rho(\mathfrak{e}_+) v_k=v_{k+1},\\
&\rho(\mathfrak{h}) v_k=(\omega+2k) v_{k}
\end{align*}
Since $v_k=\mathfrak{e}_+^k (v_0)$, $W_\omega$ is completely determined by its lowest weight vector $v_0$, i.e., the unique generating element of the kernel of the operator $\rho(\mathfrak{e}_-)$ in $W_\omega$. The first step in the algorithm for decomposing an arbitrary $\mathfrak{sl}_2$-module into irreducible components is finding all its linearly independent lowest weight vectors, see \cite{FH}.

Consider the following representation by differential operators on $\mathbb{Q}[x_1, \ldots, x_n]$:
\begin{align*}
&D_+=\rho_1(\mathfrak{e}_+)=x_1^2\partial_{1}+x_2^2\partial_{2}+\cdots+x_n^2\partial_{n},\\
&D_0=\rho_1(\mathfrak{h})=2(x_1\partial_{1}+x_2\partial_{2}+\cdots+x_n\partial_{n}),\\
&D_-=\rho_1(\mathfrak{e}_-)=-(\partial_{1}+\partial_{2}+\cdots+\partial_{n}),
\end{align*}
where $\partial_i=\dfrac{\partial}{\partial x_i}$.

Direct verification shows that these operators satisfy the following commutation relations
$$
[D_+,D_-]=D_0, [D_0, D_+]=2 D_+, [D_0,D_-]=-2D_-,
$$
and indeed define an $\mathfrak{sl}_2$-action on the algebra $\mathbb{Q}[x_1, \ldots, x_n]$.
Each of the operators maps a symmetric polynomial back to a symmetric polynomial, so the vector space $\Lambda_n$ is an $\mathfrak{sl}_2$-invariant subspace of the representation $\rho_1$. From now on, we consider the operators $D_\pm, D_0$ only on $\Lambda_n$.

To find the decomposition of $\Lambda_n$ into irreducibles, we need to describe all lowest weight vectors in $\Lambda_n$, i.e., find the kernel $\ker D_-=\{ z \in \Lambda_n \colon D_-(z)=0 \}$ of the operator $D_-$. This operator is a locally nilpotent derivation on $\Lambda_n$, and for such operators, the algorithm for finding the kernel is well-known. Let $D$ be a locally nilpotent derivation of a commutative $\mathbb{Q}$-algebra $A$, generated by elements $a_1, a_2, \ldots, a_n$. Suppose that $D$ has a \textit{slice}, i.e., an element $s \in A$, for which $D(s)=-1$ holds.
Then $\ker D$ is generated by elements $\sigma(a_1), \sigma(a_1), \ldots $ where the homomorphism $\sigma: A \longrightarrow \ker D$ is defined as
$$
\sigma(a)=\sum_{i=0}^\infty \frac{1}{i!}  D^i(a) s^i,
$$
see \cite{Essen} for details.

In the following theorem, using this technique of locally nilpotent derivations, we explicitly describe the kernels of the derivations $D_-, D_+$, and also present a new basis of $\Lambda_n$ formed from the lowest weight vectors of the representation $\rho_1$.
\begin{te}\label{dp}  The following holds
\begin{enumerate}[label={\rm (\roman*)}, left=0pt]
\item $
\ker D_-=\mathbb{Q}[z_2,z_3,\ldots,z_n],
$
where $$
z_i= \sum_{k=0}^{i-2} (-1)^k n^{i-k-1} {i \choose k} p_{i-k} p_1^k  +(i-1)(-1)^{i+1} p_1^i, i=2,3,\ldots, n,
$$
and $p_k=x_1^k+x_2^k+\cdots+x_n^k$ is the power sum symmetric polynomial.
\item  $
\ker D_+=\mathbb{Q}.
$
\item $\Lambda_n=\mathbb{Q}[p_1, z_2,z_3,\ldots,z_n].$
\end{enumerate}
\end{te}
\begin{proof}
\textit{ (i)} Let's choose the  basis in the vector space $\Lambda_n$ of  the power sum symmetric polynomials $p_i=x_1^i+x_2^i+\cdots+x_n^i, i=1,2,\ldots, n$. The operator $D_-$ is locally nilpotent on $\Lambda_n$ with the following action on the basis elements: $D_-(p_i)=-i p_{i-1}, D_-(p_1)=-n.$ Since
$$
D_-\left(\frac{p_1}{n}\right)=-1,
$$
we have a slice.
Then the following elements
$$
\sigma(p_i)=\sum_{k=0}^\infty \frac{1}{k!}D_-^k(p_i) \left(\frac{p_1}{n} \right)^k,  2 \leq i \leq n,
$$
generate the kernel $\ker D_-$, see \cite[Corolary 1.3]{Essen1}. Since the operator $D_-$ is locally nilpotent, the sum is finite. After some simplifications, we get an explicit expression for the generating elements of the kernel:
$$
z_i= \sum_{k=0}^{i-2} (-1)^k n^{i-k-1} {i \choose k} p_{i-k} p_1^k  +(i-1)(-1)^{i+1} p_1^i, i=2,3\ldots,n.
$$
Let  us  prove that the elements $z_i$ are not zero for $n>2$. Denote by $C_i$ the coefficient of $x_1^i$ in $z_i$, which obviously equals
$$
C_i=\sum_{k=0}^{i-2} (-1)^k {i \choose k} n^{i-k-1} +(i-1)(-1)^{i+1}.
$$
After simple transformations, we get
$$
n C_i=(n-1)^i+(-1)^i(n-1).
$$
It's easy to show that the equation $(n-1)^i+(-1)^i(n-1)=0$ has only two integer solutions $n=0,n=1$ for even $i$ and three integer solutions $n=0,n=1,n=2$ for odd $i$.
Therefore, the coefficient of $x_1^i$ in $z_i$ for $n>2$ is non-zero, so and $z_i$ is non-zero.

\textit{(ii)}  The operator $D_+$ acts on power sums as follows: $D_+(p_i)=-i p_{i+1}$. Let $\tilde p_i =-(i-1)! p_i.$ Then
$D_+(\tilde p_i)=\tilde p_{i+1}.$  From \cite[Proposition 3.2.5.]{Now} it follows that such a derivation has a trivial kernel, which is what we needed to prove.

\textit{(iii)} The algebra $\Lambda_n$ is a polynomial ring in the  slice  over $\ker D_-$, see \cite[Proposition 2.1]{Wr}. Therefore, taking into account $(i)$, we get that
$$ \Lambda_n=\mathbb{Q}[p_1, z_2,z_3,\ldots,z_n]. $$
We have obtained another family of symmetric polynomials which forms a basis in $\Lambda_n$.
\end{proof}
Let's provide explicit expressions for $z_i$:
\begin{pr}{\rm
\begin{align*}
&z_2=np_{{2}}-{p_{{1}}}^{2},\\
&z_3={n}^{2}p_{{3}}-3\,np_{{1}}p_{{2}}+2\,{p_{{1}}}^{3},\\
&z_4={n}^{3}p_{{4}}-4\,{n}^{2}p_{{1}}p_{{3}}+6\,n{p_{{1}}}^{2}p_{{2}}-3\,{p
_{{1}}}^{4},\\
&z_5={n}^{4}p_{{5}}-5\,{n}^{3}p_{{1}}p_{{4}}+10\,{n}^{2}{p_{{1}}}^{2}p_{{3}
}-10\,n{p_{{1}}}^{3}p_{{2}}+4\,{p_{{1}}}^{5}.
\end{align*}
}
\end{pr}
As a vector space, the kernel $\ker D_-$ is generated by monomials $z^\alpha=z_2^{\alpha_1} z_3^{\alpha_2} \cdots z_n^{\alpha_{n-1}} $ $\alpha \in \mathbb{N}^{n-1}$. Denote by $W(z^\alpha)$ the vector space generated by polynomials $D_+^k(z^\alpha)$, $k \in \mathbb{N}$.
The following theorem gives an explicit decomposition of $\Lambda_n$ into irreducible $\mathfrak{sl}_2$-modules.
\begin{te}\label{dc}
The vector space $
W(z^\alpha) ,
$
is isomorphic to the standard lowest weight  $\mathfrak{sl}_2$-module $V_{w(\alpha)}$ with the weight $$
w(\alpha)=2 (2 \alpha_1+ 3 \alpha_2+\cdots+n \alpha_{n-1}).
$$
For $n>2$ we have a direct sum decomposition
$$
\Lambda_n = \bigoplus_{\alpha \in \mathbb{N}^{n-1}}  W(z^\alpha),
$$
and 
$$
\Lambda_n \cong \bigoplus_{i=0}^\infty c_i W_{2i},
$$
where $c_i$ is the number of natural solutions to the equation $2 \alpha_1+ 3 \alpha_2+\cdots+n \alpha_{n-1}=i$ and $W_{2i}$ is the standard infinite dimensional $\mathfrak{sl}_2$-module of the type $([\circ )$ with the weight $2i$.
\end{te}
\begin{proof}
Let's show that $W(z^\alpha)$ is isomorphic to the standard lowest weight infinite-dimesional module $V_\omega$, for some weight $\omega$.
Consider in $W(z^\alpha)$ a basis formed by vectors $v_i=D_+(z^\alpha)$. Since the kernel of the operator $D_+$ is trivial (Theorem~\ref{dp},$(ii)$), $W(z^\alpha)$ is an  infinite-dimensional vector space. From $D_0(z_i)=2i z_i$ we get $D_0(z^\alpha)=w(\alpha) z^\alpha$, where
$$
w(\alpha)=2 (2 \alpha_2+ 3 \alpha_3+\cdots+n \alpha_n).
$$
By direct calculations in the Lie algebra $\mathfrak{sl}_2$, it can be shown that for any element $z \in \ker D_-$ of the  weight $w(z)$, the following relations hold
\begin{align*}
&D_0 D_+^i (z)=(w(z)+2i)D_+^i (z),\\
&D_- D_+^i (z)=-i(w(z)+i-1)D_+^{i-1} (z)
\end{align*}
Consider in $W(z^\alpha)$ a basis formed by vectors $u_i=D_+(z^\alpha)$. Then we have that
\begin{align*}
&D_0(v_i)=(w(\alpha)+2i) v_i,\\
&D_-(v_i)=-i(w(\alpha)+i-1) u_{i-1},\\
&D_+(v_i)=v_{i+1}
\end{align*}
Therefore, $W(z^\alpha)$ is isomorphic to the standard irreducible lowest weight module  $W_{	w(\alpha)}$.  Then 
$\Lambda_n$ is a direct sum of lowest weight modules generated by linearly independent  lower weight vectors. The vectors of lower weight in $\Lambda_n$ are linearly independent elements of the kernel of the derivation $D_+$, i.e., monomials of the form $z^\alpha=z_2^{\alpha_1} z_3^{\alpha_2} \cdots z_n^{\alpha_{n-1}} $ where $\alpha$ runs through all tuples of natural numbers of length $n-1$.
Thus, we have a decomposition into irreducibles modules
$$
\Lambda_n= \bigoplus_{\alpha \in \mathbb{N}^{n-1}}  W(z^\alpha).
$$
By virtue of the isomorphism $ W(z^\alpha)$ and $W_{w(\alpha)}$ we have
$$
\Lambda_n \cong \bigoplus_{\alpha \in \mathbb{N}^{n-1}}  W_{w(\alpha)}.
$$
When $\alpha$ runs through  $\mathbb{N}^{n-1}$, $w(\alpha)$ many  times runs through all even numbers.
It's clear that the multiplicity $c_i$ of the occurrence of $W_{2i}$ in $\Lambda_n$ is equal to the number of natural solutions to the equation $w(\alpha)=2i$. Then
$$
\Lambda_n \cong \bigoplus_{i=0}^\infty c_i W_{2i}.
$$
\end{proof}
\begin{pr}{\rm  Consider the case $n=3$. We have $\ker D_-=\mathbb{Q}[z_2,z_3]$, where

\begin{align*}
&z_2=3 p_{{2}}-{p_{{1}}}^{2}=2\,{x_{{1}}}^{2}-2\,x_{{1}}x_{{2}}-2\,x_{{1}}x_{{3}}+2\,{x_{{2}}}^{2}-2\,x_{{2}}x_{{3}}+2\,
{x_{{3}}}^{2},\\
&z_{{3}}=2\,{p_{{1}}}^{3}-9\,p_{{1}}p_{{2}}+9\,p_{{3}}=\left( x_{{3}}+x_{{1}}-2\,x_{{2}} \right)  \left( x_{{1}}+x_{{2}} -2\,x_{{3}}\right)  \left( 2\,x_{{1}}-x_{{3}}-x_{{2}} \right).
\end{align*}

Then
$$
\Lambda_3=\bigoplus_{\alpha \in \mathbb{N}^2} W(z_2^{\alpha_1} z_3^{\alpha_2}).
$$
Let's gather the isomorphic components
$$
\Lambda_3\cong\bigoplus_{2\alpha_1+3 \alpha_2=i} W_{2i}=\bigoplus_{i=0}^\infty  c_i W_{2i},
$$
where $c_i$ is the number of natural solutions to the equation $2\alpha_1+3 \alpha_2=i$. For the sequence $c_i$  its  generating function is
$$
\frac{1}{(1-x^2)(1-x^3)}=\frac{1}{{z}^{5}-{z}^{3}-{z}^{2}+1}=1+{z}^{2}+{z}^{3}+{z}^{4}+{z}^{5}+2,{z}^{6}+{z}^{7}+2\,{z}^{8}+2\,{z}^{9}+2\,{z}^{10}+\cdots.
$$
We obtain the following recurrence relation
$$
c_{i-5}-c_{i-3}-c_{i-2}+c_n=0, c_0=1,c_1=0,c_2=1,c_3=1,c_4=1.
$$
Solving which we get
$$
c_i=\frac 1 6 i+\frac 1 4 (-1)^i+{\frac {5}{
12}}+\frac 1 3\,\cos \left( \frac{2\,\pi \,i}{3}
\right) -\frac 1 9\,\sqrt {3}\sin \left(\frac{ 2\,\pi \,i}{3} \right).
$$
}
\end{pr}
Let  us  find the decomposition for $n=2$. The kernel $D_-$ in $\Lambda_2$ can be found directly, it is generated by the powers of the symmetric polynomial $z_2=(x_1-x_2)^2$. For all powers $i$ we have
$$
D_+(z_2^{2i})=2i (x_1+x_2)z_2^{2i} \neq 0.
$$
Hence, the element $z_2^{2i}, i>0$ generates an infinite-dimensional $\mathfrak{sl}_2$-module $ W(z_2^{2i})$. Since $\omega(z_2^{2i})=8i$, $ W(z_2^{2i})$ is isomorphic to the standard $\mathfrak{sl}_2$-module $W_{8i}$.
Therefore, we obtain the following decomposition
$$
\Lambda_2= \bigoplus_{i=0}^\infty  W(z_2^{2i})\cong \bigoplus_{i=0}^\infty  W_{8i}.
$$
Note that for odd $i$ we get that $z_i$ is equal to zero, so we have the following identities
$$
\sum_{k=0}^{i-2} (-1)^k {i \choose k} p_{i-k} p_1^k 2^{i-k-1} =(i-1)(-1)^{i} p_1^i,
$$
or, after simplification
$$
\sum_{k=0}^{2m-1} \left(-\frac{1}{2}\right)^{k+1} {2m+1 \choose k} (x^{2m+1-k}+y^{2m+1-k}) (x+y)^k  =\frac{m}{2^{2m}}(x+y)^{2m+1}, m=1,2,\ldots.
$$
For even $i$, $z_i$ is a homogeneous polynomial of degree $2i$ in $(x_1-x_2)^i$. Since for $n=2$ the coefficient $C_{2i}$ of $x_1^{2i}$ is equal to 1, we get the following identity
$$
\sum_{k=0}^{2m-2} (-1)^k 2^{2m-k-1} {2m\choose k} (x^{2m-k}+y^{2m-k}) (x+y)^k  =(2m-1) (x+y)^{2m}+(x-y)^{2m}.
$$

\section{Action on the Schur polynomials}


Recall that the symmetric Schur polynomials $\boldsymbol{s}_\lambda$ for each partition $\lambda \in \mathcal{Y}_n$ are defined as the ratio of determinants:
$$
\boldsymbol{s}_\lambda=\frac{a_{\lambda+\delta}}{a_{\delta}},
$$ 
where 
$$
a_\mu=\det(x_i^{\mu_j})_{1\leq i,j \leq n}= \begin{vmatrix} x_1^{\mu_1} & x_2^{\mu_1} & \ldots & x_n^{\mu_1}\\
x_1^{\mu_2} & x_2^{\mu_2} & \ldots & x_n^{\mu_2}\\
\hdotsfor{4}\\
x_1^{\mu_n} & x_2^{\mu_n} & \ldots & x_n^{\mu_n}
\end{vmatrix}, 
$$
and   $\delta=(n-1,n-2, \ldots, 1,0)$.
Since the Schur polynomials 
form a basis in $\Lambda_n$ (as an infinite-dimensional $\mathbb{Q}$-vector space), the $\mathfrak{sl}_2$-action on Schur polynomials is a linear combination of Schur polynomials.

In the following theorem, the action of the operators $D_-, D_0, D_+$ on Schur polynomials is expressed in this basis. Let $c(\square)=j-i$ be the content of the cell $\square=(i,j).$

\begin{te}\label{mt} The following relations hold:
$$
\begin{array}{ll}
 (i)  &D_-(\boldsymbol{s}_\lambda)=-\sum\limits_{\mu=\lambda-\square \in\mathcal{Y}_n } (n+c(\square))\boldsymbol{s}_\mu,\\
 (ii) &D_0(\boldsymbol{s}_\lambda)=2 |\lambda| \boldsymbol{s}_\lambda, \\
 (iii)  &D_+(\boldsymbol{s}_\lambda)=\sum\limits_{\mu=\lambda+\square  \in \mathcal{Y}_n} c(\square) \boldsymbol{s}_\mu.
\end{array}
$$
\end{te}
 \begin{proof} $(i)$ Since the operator $D_-$ is a derivation of the algebra of polynomials, we have 
$$
D_-(\boldsymbol{s}_\lambda)=D_-\left(\frac{a_{\lambda+\delta}}{a_{\delta}}\right)=\frac{D_-(a_{\lambda+\delta})a_{\delta}-a_{\lambda+\delta}D_-(a_{\delta})}{a_{\delta}^2}.
$$
According to the determinant derivative  rule, we find the action of $D_-$ on the Vandermonde determinant:
\begin{gather*}
D_-(a_{\delta})=\\=\!\begin{vmatrix} D_-(x_1^{n-1}) & D_-(x_2^{n-1}) & \ldots & D_-(x_n^{n-1}) \\  x_1^{n-2} & x_2^{n-2} & \ldots & x_n^{n-2} \\ 
\vdots & \vdots & \ldots & \vdots \\
1 & 1 & \ldots & 1
\end{vmatrix}{+}\begin{vmatrix} x_1^{n-1} & x_2^{n-1} & \ldots & x_n^{n-1} \\  D_-(x_1^{n-2}) & D_-(x_2^{n-2}) & \ldots & D_-(x_n^{n-2}) \\ 
\vdots & \vdots & \ldots & \vdots \\
1 & 1 & \ldots & 1
\end{vmatrix}+\\+\cdots+\begin{vmatrix} x_1^{n-1} & x_2^{n-1} & \ldots & x_n^{n-1} \\  x_1^{n-2} & x_2^{n-2} & \ldots & x_n^{n-2} \\ 
\vdots & \vdots & \ldots & \vdots \\
0 & 0 & \ldots & 0
\end{vmatrix}.
\end{gather*}
Since $D_-(x_i^k)=-\partial_i(x_i^k)=-k\, x_i^{k-1}$, in each of these determinants, except for the last one, there are two proportional rows, and in the last determinant, there is a zero row. Therefore, the derivarion  $D_-$ annihilates the Vandermonde determinant, and we obtain a simpler expression for the action of the operator $D_-$:
$$
D_-(\boldsymbol{s}_\lambda)=\frac{D_-\left(\begin{vmatrix} x_1^{\lambda_1+n-1} & x_2^{\lambda_1+n-1} & \ldots & x_n^{\lambda_1+n-1} \\  x_1^{\lambda_2+n-2} & x_2^{\lambda_2+n-2} & \ldots & x_n^{\lambda_2+n-2} \\ 
\vdots & \vdots & \ldots & \vdots \\
x_1^{\lambda_n} & x_2^{\lambda_n} & \ldots & x_n^{\lambda_n}
\end{vmatrix}\right)}{a_{\delta}}.
$$
We differentiate the determinant in the numerator and obtain  the following sum
\begin{gather*}
D_-\left(a_{\lambda+\delta}\right)=-\sum_{i=1}^n(n{+}\lambda_i{-}i)\!\begin{vmatrix} x_1^{\lambda_1+n-1} & x_2^{\lambda_1+n-1} & \ldots & x_n^{\lambda_1+n-1} \\
\vdots & \vdots & \ldots & \vdots \\
 x_1^{\lambda_i+n-i-1} & x_2^{\lambda_i+n-i-1} & \ldots & x_n^{\lambda_i+n-i-1}\\
\vdots & \vdots & \ldots & \vdots \\
x_1^{\lambda_n} & x_2^{\lambda_n} & \ldots & x_n^{\lambda_n}
\end{vmatrix}\!.
\end{gather*}

It is clear that each of the determinants in the sum is nonzero only when the cell $(i,\lambda_i)$ in the Young diagram corresponding to the partition $\lambda$ is an \textit{inner cell}. Otherwise, we get a determinant with two proportional rows. Note that $\lambda_i-i$ is the content $c(\square)$ of this inner cell. Therefore, dividing this sum by the Vandermonde determinant, we obtain  
$$
D_-(\boldsymbol{s}_\lambda)=-\sum_{\mu=\lambda-\square \in\mathcal{Y}_n } (n+c(\square)) s_{\mu}.
$$

Note that the action of the operator $D_-$ on the Schur polynomials was recently considered in \cite{W23}, and in \cite{Gr}, the action of $D_-$ was generalized to the skew Schur polynomials.

$(ii)$ We now find the action of the derivation  $D_+$ on the Vandermonde determinant. Since $D_+(x_i^k)=k x_i^{k+1}$, in each determinant of the sum $D_+(a_{\delta})$, except for the one corresponding to the differentiation of the first row, there will be proportional rows. As a result, only one determinant remains:
\begin{gather*}
D_+\left(\begin{vmatrix} x_1^{n-1} & x_2^{n-1} & \ldots & x_n^{n-1} \\  x_1^{n-2} & x_2^{n-2} & \ldots & x_n^{n-2} \\ 
\vdots & \vdots & \ldots & \vdots \\
1 & 1 & \ldots & 1
\end{vmatrix} \right)=(n-1) \begin{vmatrix} x_1^{n} & x_2^{n} & \ldots & x_n^{n} \\  x_1^{n-2} & x_2^{n-2} & \ldots & x_n^{n-2} \\ 
\vdots & \vdots & \ldots & \vdots \\
1 & 1 & \ldots & 1
\end{vmatrix}.
\end{gather*}
The resulting determinant is a \textit{generalized Vandermonde determinant}.
Using the formula for its computation, see \cite[Lemma 2.1]{Ernst}, we obtain
$$
\begin{vmatrix} x_1^{n} & x_2^{n} & \ldots & x_n^{n} \\  x_1^{n-2} & x_2^{n-2} & \ldots & x_n^{n-2} \\ 
\vdots & \vdots & \ldots & \vdots \\
1 & 1 & \ldots & 1
\end{vmatrix}=
(x_1+x_2+\cdots+x_n) a_{\lambda+\delta}=\boldsymbol{s}_{(1)} a_{\delta}.
$$
Applying the derivation  $D_+$ to the Schur polynomial, we obtain 
\begin{gather*}
D_+(\boldsymbol{s}_\lambda)=D_+\left(\frac{a_{\lambda+\delta}}{a_{\delta}}\right)=\frac{D_+(a_{\lambda+\delta})a_{\delta}-a_{\lambda+\delta}D_+(a_{\delta})}{a_{\delta}^2}=\frac{D_+\left(a_{\lambda+\delta} \right)}{a_{\delta}}-(n-1)\boldsymbol{s}_{(1)} \boldsymbol{s}_\lambda.
\end{gather*}

The expression $D_+\left(a_{\lambda+\delta} \right)$ is a sum of determinants that are nonzero only when the cell $(i,\lambda_i+1)$ is an \textit{outer corner} of the Young diagram $\lambda.$
We have
\begin{gather*}
\frac{D_+\left(\det(x_j^{\lambda_i+n-i}) \right)}{a_{\delta}}=\frac{\sum\limits_{j=1}^n(\lambda_i+n-i)\det(x_j^{\lambda_i+n-i})}{a_{\delta}}=\sum_{\mu=\lambda+\square \in \mathcal{Y}_n} (n-1+c(\square)) \boldsymbol{s}_\mu.
\end{gather*}
According to the Pieri rule, we find
$
\boldsymbol{s}_{(1)} \boldsymbol{s}_\lambda=\sum\limits_{\mu=\lambda+\square \in \mathcal{Y}_n} \boldsymbol{s}_\mu.
$ 

In conclusion, we obtain
$$
D_{+}(\boldsymbol{s}_{\lambda})=\sum_{\mu=\lambda+\square  \in \mathcal{Y}_n} c(\square) \boldsymbol{s}_\mu.
$$
$(iii)$ We find  the action of the operator  $D_0$ on the Vandermonde determinant
\begin{gather*}
D_0 \!\left(\begin{vmatrix} x_1^{n-1} & x_2^{n-1} & \ldots & x_n^{n-1} \\  x_1^{n-2} & x_2^{n-2} & \ldots & x_n^{n-2} \\ 
\vdots & \vdots & \ldots & \vdots \\
1 & 1 & \ldots & 1
\end{vmatrix}  \right)\! \!  =2 ((n{-}1){+}(n{-}2){+}\cdots{+}1)\begin{vmatrix} x_1^{n-1} & x_2^{n-1} & \ldots & x_n^{n-1} \\  x_1^{n-2} & x_2^{n-2} & \ldots & x_n^{n-2} \\ 
\vdots & \vdots & \ldots & \vdots \\
1 & 1 & \ldots & 1
\end{vmatrix}\!\!=\\=n(n-1) a_{\delta}.
\end{gather*}

Then 
\begin{gather*}
D_0(\boldsymbol{s}_{\lambda})=D_0 \left( 
\frac{\det(x_j^{\lambda_i+n-i})}{a_{\delta}}\right)=\frac{D(\det(x_j^{\lambda_i+n-i})) a_{\delta}-\det(x_j^{\lambda_i+n-i}) D_0(a_{\delta})}{a_{\delta}^2}=\\=
\frac{ \sum\limits_{i=1}^n 2 (\lambda_i+n-i) \det(x_j^{\lambda_i+n-i}) a_{\delta}-\det(x_j^{\lambda_i+n-i}) n(n-1) a_{\delta}}{a_{\delta}^2}=2 |\lambda| \boldsymbol{s}_\lambda.
\end{gather*}
\end{proof}

As an immediate consequence of Theorem~\ref{mt}, we obtain the explicit form of the action of the operators $D_-, D_0, D_+$ on various bases in $\Lambda_n$ 

\begin{te}\label{t5} The action of the operators 
\begin{enumerate}[label={\rm \arabic*)}] 
\item
on the elementary symmetric polynomials $e_i$ and on the corresponding partitions $(1^i):$

 \begin{align*}
&D_-(e_i)=  -(n-(i-1))e_{i-1}\\
&D_0(e_i)=2i e_i,\\
&D_+(e_i)=e_1 e_i - (i+1) e_{i+1}, \\
&D_+(e_n)=e_1 e_n.
\end{align*}

\item   On the complete  homogeneous symmetric polynomials $h_i$:
 \begin{align*}
 &D_-(h_i)=-(n+i-1)h_{i-1}, \\
 &D_0(h_i)=2i h_i,\\
&D_+(h_i)=(i+1) h_{i+1}-h_1 h_i.
\end{align*}

\item  On the power symmetric polynomials $p_i$:
 \begin{align*}
&D_-(p_i)=i p_{i-1},\\
&D_0(p_i)=2i p_i,\\
&D_+(p_i)=i p_{i+1}, 
\end{align*}
\end{enumerate}
for all $1 \leq i \leq n.$
\end{te} 
\begin{proof}
  From Theorem~\ref{mt}, since the symmetric polynomials $e_i$ and $h_i$ are Schur polynomials corresponding to the partitions $(1^i)$ and $(i)$, we immediately obtain that 
$$
\begin{array}{ll}
D_-(e_i)=  -(n-(i-1))e_{i-1} & D_-(h_i)=-(n+i-1)h_{i-1}\\
D_0(e_i)=2i e_i & D(h_i)=2i h_i,\\
D_+(e_i)=e_1 e_i - (i+1) e_{i+1} & D_+(h_i)=(i+1) h_{i+1}-h_1 h_i,\\
D_+(e_n)=e_1 e_n & D_+(h_n)=(n+1) h_{n+1}-h_1 h_n.\\
\end{array}
$$
The results of the action of the operators $D_-, D_0$ are clear, let us explain the obtained expressions for $D_+$.
By Theorem~\ref{mt}, we have that 
$$ D_+(\boldsymbol{s}_\lambda)=\sum\limits_{\mu=\lambda+\square  \in \mathcal{Y}_n} c(\square) \boldsymbol{s}_\mu.$$
For $\lambda=(1^i)$, $i<n$, we obtain, 
\begin{gather*}
D_+(\boldsymbol{s}_{(1^i)})=D_+(e_i)=\sum\limits_{\mu=(1^i)+\square  \in \mathcal{Y}_n} c(\square) \boldsymbol{s}_\mu=\boldsymbol{s}_{2, 1,\ldots,1}- i e_{i+1}.
\end{gather*}
According to the Pieri rule, we have 
$$ e_1 e_i=\boldsymbol{s}_{2, 1,\ldots,1}+e_{i+1}.$$
Hence, 
$$
D_+(e_i)=e_1 e_i - (i+1) e_{i+1}.
$$
For the case $i=n$, we have $
D_+(e_n)=e_1 e_n.
$

For the complete homogeneous symmetric polynomials, the expression for $D_+(h_i)$ is obtained similarly.
The action on the power symmetric polynomials $p_i$ is immediately obtained from the definition of these polynomials.

\end{proof}

\section{The second $\mathfrak{sl}_2$-action on $\Lambda_n$. }


Consider the  $\mathbb{Q}$-vector space $V_d$ of dimension $d+1$ with the  basis $\langle  v_0, v_1, \ldots, v_d \rangle.$ The following linear operators 
\begin{align*}
&\rho(\mathfrak{e}_+)(v_i)=(d-i) v_{i+1},\\
&\rho(\mathfrak{h})(v_i)=(2i-d) v_i,\\
&\rho(\mathfrak{e}_-)(v_i)=i v_{i-1},
\end{align*}
turn $V_d$ into an irreducible finite-dimensional $\mathfrak{sl}_2$-module, which is called the finite-dimensional \textit{standard} $\mathfrak{sl}_2$-module. In the terminology of \cite{HT}, this is a module of type $([\circ])$.

For any $\mathfrak{sl}_2$-module $V$, the eigenvalues of the operator $\rho(\mathfrak{h})$ are called the \textit{weights} of $V$. The weight of the element $v_0$, i.e., the unique basis element of the kernel $\ker \rho(\mathfrak{e}_-)$ is equal to $\dim V_d-1.$ Any $\mathfrak{sl}_2$-module, up to isomorphism, is determined by the \textit{multiset} of its weights. The formal sum 
$$
{\rm Char} V=\sum_{w} q^w,
$$
where $w$ runs over the multiset of weights of $V$, is called the \textit{character} of the $\mathfrak{sl}_2$-module $V$. The character determines any finite-dimensional $\mathfrak{sl}_2$-module up to isomorphism, see \cite{FH}.

Let $\mathcal{P}_n^{(d)}$ be the subspace in the vector space $\mathbb{Q}[x_1,x_2, \ldots, x_n]$ generated by polynomials whose degree in each variable does not exceed $d$. The basis of $\mathcal{P}_n^{(d)}$ consists of monomials of the form $x_1^{k_1} x_2^{k_2} \cdots x_n^{k_n}$, with $k_i \leq d$, so its dimension is $(1+d)^n$.
 
The symmetric group $S_n$ acts by permutations on $\mathcal{P}_n^{(d)}$, and let $\Lambda_n^{(d)}$ be the vector space of all symmetric polynomials in $\mathcal{P}_n^{(d)}$. It is clear that the basis of $\Lambda_n^{(d)}$ consists of all Schur polynomials $\boldsymbol{s}_\lambda$ for which $\lambda \leq (\underbrace{d,d,\ldots,d}_{n \text{ times}}).$

Let us write the partition $\lambda=(\lambda_1, \ldots, \lambda_n), \lambda_1 \leq d$ in the notation $\lambda=(0^{\alpha_0} 1^{\alpha_1} \cdots d^{\alpha_d}), |\alpha|=n.$ Since the complete homogeneous symmetric polynomial $h_n(y_0,y_1,\ldots,y_d)$ is the sum of monomials of the form 
$ y_0^{\alpha_0} y_1^{\alpha_1} \cdots y_d^{\alpha_d},|\alpha|=n, 
$
the dimension of $\Lambda_n^{(d)}$ is equal to the specialization 
$$
\dim \Lambda_n^{(d)} =h_{n}(\underbrace{1,1,\ldots,1}_{d+1})=\binom{n+d}{n}.
$$

It can be easily demonstrated that the following operators on $\mathcal{P}_n^{(d)}$:
\begin{align*}
&D_+ =\rho_2(\mathfrak{e}_+)=\sum_{k=1}^n(- x_k^2\partial_{k}+d \cdot x_k),\\ 
&D_0 = \rho_2(\mathfrak{h})=2 \sum_{k=1}^n x_k\partial_{k}-n \cdot d, \\
&D_- =\rho_2(\mathfrak{e}_-) =\sum_{k=1}^n \partial_{k}.
\end{align*}
satisfy the commutation relations 
$$
[D_+,D_-]=D_0, [D_0, D_+]=2 D_+, [D_0,D_-]=-2D_-.
$$
and therefore define an $\mathfrak{sl}_2$-module structure on $\mathcal{P}_n^{(d)}$, see also \cite{KO}.

Moreover, they map symmetric polynomials to symmetric ones, so $\Lambda_n^{(d)}$ is also an $\mathfrak{sl}_2$-module.

The following theorem establishes the structure of the module.

\begin{te}\label{fc}  The following decomposition holds:
$$
\Lambda_n^{(d)} \cong \bigoplus\limits_{i=0}^{nd} c_n(d,i) V_i,
$$
 where  $$
c_n(d,i)=\gamma \left( d + n, n, \frac{1}{2}dn - \frac 1 2 i \right) - \gamma \left( d + n, n, \frac{1}{2}dn - \frac 1 2 i - 1 \right).
$$
and $\gamma(d,n,i)$ is  the number of
partitions of size $i$  inside the rectangle $d \times  n$.
\end{te} 

\begin{proof}  First, we will prove that $\Lambda_n^{(d)} \cong {\rm Sym}^n (V_d)$, where ${\rm Sym}^n (V_d)$ is the  $n$-graded component of the \textit{symmetric algebra} ${\rm Sym} (V_d)$. To prove the isomorphism of $\mathfrak{sl}_2$-modules $\Lambda_n^{(d)}$ and ${\rm Sym}^n (V_d)$, it is sufficient to show that they have the same characters. 
We choose the standard basis $V_d=\langle v_0, v_1, \ldots, v_d \rangle$ in which $D_0(v_i)=2i-d.$ The symmetric power ${\rm Sym}^n (V_d)$ is generated by the monomials of the complete symmetric polynomial $h_n(v_0, v_1, \ldots, v_d )$. Since $\rho(\mathfrak{h})$ is a derivation  on ${\rm Sym}^n (V_d)$, for each such monomial we have 
\begin{gather*}
\rho(\mathfrak{h}) (v_0^{\alpha_0}, v_1^{{\alpha_1}}, \ldots, v_d^{\alpha_d})=\left(\sum_{i=0}^d \alpha_i(2i-d)\right) v_0^{\alpha_0}  v_1^{{\alpha_1}}  \cdots, v_d^{\alpha_d}.
\end{gather*}
Therefore, 
$$
{\rm Char}{\rm Sym}^n (V_d)=\sum_{|\alpha|=n } q^{2 \sum\limits_{i=0}^d i \alpha_i - nd}.
$$
Recall that the basis of $\Lambda_n^{(d)}$ consists of the Schur polynomials $\boldsymbol{s}_{\lambda}$, where $\lambda$ runs over all partitions of length $n$, including zero parts, with $\lambda_1 \leq d.$ We write $\lambda$ in the notation $\lambda=(0^{\alpha_0} 1^{\alpha_1} \cdots d^{\alpha_d}).$ Then 
$$
|\lambda|=\sum_{i=0}^d i \alpha_i, |\alpha|=n.
$$  
The eigenvalue of $D_0$ on $\boldsymbol{s}_{\lambda}$, which equals $2 |\lambda|-nd$, coincides with the eigenvalue of $\rho(\mathfrak{h})$ on $v_0^{\alpha_0}, v_1^{{\alpha_1}}, \ldots, v_d^{\alpha_d}.$ 
Therefore, 
$$
{\rm Char}\left( \Lambda_n^{(d)}\right)=\sum_{\lambda \leq (d,d, \ldots,d)} q^{2 |\lambda|-nd}=\sum_{|\alpha|=n } q^{2 \sum\limits_{i=0}^d i \alpha_i - nd}={\rm Char}\left( {\rm Sym}^n (V_d) \right).
$$
Hence $\Lambda_n^{(d)} \cong {\rm Sym}^n (V_d).$

The problem of decomposing ${\rm Sym}^n (V_d)$ into irreducibles was one of the important tasks in classical invariant theory.  It is a classical
result  that the decomposition  into irreducible
$\mathfrak{sl}_2$-modules can be computed in terms of partitions. We will provide only the final result. 

Let $\gamma(d,n,i)$ denotes the number of
partitions of size $i$  inside the rectangle $d \times  n$,  i.e.  the coefficient of $T^i$ in the following expression:
$$
\frac{(1-T^{d-n+1}) \cdots (1-T^d)}{(1-T)\cdots (1-T^n)}.
$$
Then 
$$
c_n(d,i)=\gamma \left( d + n, n, \frac{1}{2}dn - \frac 1 2 i \right) - \gamma \left( d + n, n, \frac{1}{2}dn - \frac 1 2 i - 1 \right).
$$
In terms of classical invariant theory, the number $c_n(d,i)$ is  the dimension of the vector space of covariants of degree $d$ and order $i$  for a binary $n$-form.  The explicit formula we have provided is known as the \textit{Sylvester-Cayley formula}; for a modern exposition, see, for example, \cite[3.3.6]{Sp}.
 
\end{proof}

\begin{pr}{\rm Consider the case $n=2, d=2$. The basis of the vector space $\Lambda_2^{(2)}$ consists of the following 6 symmetric polynomials:
 $$  \boldsymbol{s}_{0,0}, \boldsymbol{s}_{1,0}, \boldsymbol{s}_{1,1}, \boldsymbol{s}_{2,0}, \boldsymbol{s}_{2,1}, \boldsymbol{s}_{2,2}. $$ 
 We find the eigenvalues of the operator $D_0$ on the basis elements:
$$
-4,-2,0,0,2,4.
$$
Thus, the character has the form 
$$
{\rm Char}\left( {\rm Sym}^2 (V_2) \right)=q^{-4}+q^{-2}+2+ q^2+q^4={\rm Char}(V_2)+{\rm Char}(V_4).
$$
The vector space ${\rm Sym}^2 (V_2)$, where $V_2=\langle v_0, v_1,v_2 \rangle$, has the following basis:
$$
v_0^2,v_0 v_1, v_1^2,v_0v_2, v_1 v_2,v_2^2.
$$
The operator $\rho(\mathfrak{h})$ acts on the basis vectors as follows: $\rho(\mathfrak{h})(v_i)=2i-2.$ As a result, we obtain the same set of eigenvalues:
\begin{align*}
&\rho(\mathfrak{h})(v_0^2)=-4,\rho(\mathfrak{h})(v_0 v_1)=-2, \rho(\mathfrak{h})(v_1^2)=0,\\
&\rho(\mathfrak{h})(v_0v_2)=0, \rho(\mathfrak{h})(v_1 v_2)=2,\rho(\mathfrak{h})(v_2^2)=4.
\end{align*}
Therefore, $\Lambda_2^{(2)} \cong {\rm Sym}^2 (V_2)\cong V_2\oplus V_4.$
}
\end{pr}

In this way, we can decompose $\Lambda_n^{(d)}$ into irreducibles:

\begin{pr}
\begin{align*}
&\Lambda_3^{(2)}= {\rm Sym}^3 (V_2)=V_6 + V_2,\\
&\Lambda_3^{(3)}= {\rm Sym}^3 (V_3)=V_{{3}}+V_{{5}}+V_{{9}},\\
&\Lambda_3^{(4)}=V_{{0}}+V_{{4}}+V_{{6}}+V_{{8}}+V_{{12}},\\
&\Lambda_3^{(5)}=V_{{3}}+V_{{5}}+V_{{7}}+V_{{9}}+V_{{11}}+V_{{15}},\\
&\Lambda_3^{(6)}=V_{{2}}+2\,V_{{6}}+V_{{8}}+V_{{10}}+V_{{12}}+V_{{14}}+V_{{18}},\\
&\Lambda_3^{(7)}=V_{{3}}+V_{{5}}+V_{{7}}+2\,V_{{9}}+V_{{11}}+V_{{13}}+V_{{15}}+V_{{17}}
+V_{{21}},\\
&\Lambda_3^{(8)}=V_{{0}}+V_{{4}}+V_{{6}}+2\,V_{{8}}+V_{{10}}+2\,V_{{12}}+V_{{14}}+V_{{
16}}+V_{{18}}+V_{{20}}+V_{{24}}
\end{align*}
\end{pr}

The theorem establishes a connection between symmetric polynomials and classical invariant theory.

We can explicitly describe the irreducible submodules of $\Lambda_n^{(d)}$.
\begin{pr}{\rm  Consider the case $n=3, d=6.$ We have 
\begin{align*}
&z_2=3 p_{{2}}-{p_{{1}}}^{2}=2\,{x_{{1}}}^{2}-2\,x_{{1}}x_{{2}}-2\,x_{{1}}x_{{3}}+2\,{x_{{2}}}^{2}-2\,x_{{2}}x_{{3}}+2\,
{x_{{3}}}^{2},\\
&z_{{3}}=2\,{p_{{1}}}^{3}-9\,p_{{1}}p_{{2}}+9\,p_{{3}}=\left( x_{{3}}+x_{{1}}-2\,x_{{2}} \right)  \left( x_{{1}}+x_{{2}} -2\,x_{{3}}\right)  \left( 2\,x_{{1}}-x_{{3}}-x_{{2}} \right).
\end{align*}

The kernel $\ker D_-$ as a vector subspace in $\Lambda_3^{(6)}$ is generated by the following 8 elements:
$$
1, z_2, z_3,z_2^2,z_{{3}}z_{{2}},{z_{{2}}}^{3}, z_3^2,{z_{{2}}}^{4}-2\,{z_{{3}}}^{2}z_{{2}}
$$
with the corresponding weights:

\begin{gather*}
w(1)=-18, w(z_2)=-14, w(z_3)=-12,  w(z_{{3}}z_{{2}})=-8, \\w(z_2^2)=-10,
w({z_{{2}}}^{3})=-6, w(z_3^2)=-6,w({z_{{2}}}^{4}-2\,{z_{{3}}}^{2}z_{{2}})=-2.
\end{gather*}

Thus,
\begin{gather*}
\Lambda_3^{(6)}=V(1)+V(z_2)+V(z_3)+V(z_2^2)+V(z_{{3}}z_{{2}})+V({z_{{2}}}^{3})+V(z_3^2)+V({z_{{2}}}^{4}-2\,{z_{{3}}}^{2}z_{{2}})\cong\\\cong
V_{{2}}+2\,V_{{6}}+V_{{8}}+V_{{10}}+V_{{12}}+V_{{14}}+V_{{18}}
\end{gather*}

}
\end{pr}

Now let us  describe a realization of the standard irreducible finite-dimensional $\mathfrak{sl}_2$-module $V_d$ by symmetric polynomials. 
\begin{te}\label{res} The standard finite-dimensional $\mathfrak{sl}_2$-module $V_d$ is realized in  $\Lambda_{d}^{(1)}$ with the following basis 
$$
w_i=\frac{i!}{(d-(i-1))(d-(i-2)) \cdots d} \, e_i, w_0=e_0=1, i=0, \ldots, d,
$$
where $e_i$ are the \textit{elementary symmetric polynomials} in $d$ variables.
\end{te}
\begin{proof}
We have $\Lambda_{d}^{(1)}\cong {\rm Sym}^1(V_d)=V_d  $.
The vector space $\Lambda_{d}^{(1)}$ is generated by the Schur polynomials $\boldsymbol{s}_\lambda$ where all $\lambda_i \leq 1.$ Hence, $\lambda$ runs over all partitions $(1^i), i \leq d$, including the trivial partition. Thus, the vector space  $\Lambda_{d}^{(1)}$ is generated by  the elementary symmetric polynomials $e_i$, $i \leq d.$
We choose the following basis:
$$
w_i=\frac{i!}{(d-(i-1))(d-(i-2)) \cdots d}\,  e_i, w_0=1, i=1,2, \ldots, d.
$$

From Theorem~\ref{t5}, it follows that the operators act on the elementary symmetric polynomials as follows:
\begin{gather*}
D_-(e_i)=(d-(i-1)) e_{i-1}, D_-(e_1)=d,\\
D_+(e_i)=(i+1) e_{i+1}, D_+(e_d)=0
\end{gather*}

We will prove that $D_-(w_i)=i w_{i-1}$ and $D_+(w_i)=(d-i) w_{i+1}, D_+(w_d)=0.$
In fact
\begin{gather*}
D_-(w_i)=\frac{i!}{(d-(i-1))(d-(i-2)) \cdots d}\,(d-(i-1)) e_{i-1}=i \frac{(i-1)!}{(d-(i-2)) \cdots d} e_{i-1}=i w_{i-1},\\
D_+(w_i)=\frac{i!}{(d-(i-1))(d-(i-2)) \cdots d} (i+1) e_{i+1}=\\=(d-i) \frac{(i+1)!}{(d-i)(d-(i-1))(d-(i-2)) \cdots d} \, e_{i+1}=(d-i) w_{i+1}.
\end{gather*}
Thus, we have obtained the realization of the standard finite-dimensional $\mathfrak{sl}_2$-module $V_d$.
\end{proof}

\section{$\mathfrak{sl}_2$-actions  on the Young diagrams}

Let $\mathcal{Y}$ be the set of all partitions of arbitrary length. Consider the infinite-dimensional vector space $\mathbb{Q}\mathcal{Y}$, which consists of formal finite sums of elements of $\mathcal{Y}$ with rational coefficients. On $\mathbb{Q}\mathcal{Y}$, we can introduce the structure of a two-parameter infinite-dimensional representation of the complex Lie algebra $\mathfrak{sl}_2$, defined by the following Kerov operators:
\begin{alignat*}{2}
U\, \lambda &= \sum_{ \mu=\lambda+\square \in \mathcal{Y}} &(z+c(\square)) \, &\mu,\\
L\, \lambda &= &(z z'+2|\lambda|)\, &\lambda,\\
D\, \lambda &= \sum_{\mu=\lambda-\square \in  \mathcal{Y}} &(z'+c(\square)) \, &\mu, 
\end{alignat*}
where $z, z' \in \mathbb{C}$.
More about the Kerov operators and their applications can be found in \cite{Pet}. Similar operators on partially ordered sets of differential operators were considered in \cite{St88}.

Consider two subsets in $\mathcal{Y}$: the set of Young diagrams with no more than $n$ rows $\mathcal{Y}_n$ and its finite subset $\mathcal{Y}^{(d)}_n$ of diagrams where each row has no more than $d$ boxes. 
Let $\mathbb{Q}\mathcal{Y}_n$ and $\mathbb{Q}\mathcal{Y}^{(d)}_n$ denote the corresponding subspaces in $\mathbb{Q}\mathcal{Y}$.
It is easy to see that the operator $U$ increases the number of rows in a diagram, so the operators $U, L, D$ no longer define an $\mathfrak{sl}_2$-action on $\mathcal{Y}_n$ and on $\mathbb{Q}\mathcal{Y}^{(d)}_n$.

We aim to endow $\mathbb{Q}\mathcal{Y}_n$ and $\mathbb{Q}\mathcal{Y}^{(d)}_n$ with the structure of $\mathfrak{sl}_2$-modules by using the already known $\mathfrak{sl}_2$-module structure on $\Lambda_n$ and $\Lambda_n^{(d)}$.  If   we define the multiplication operation on $\mathbb{Q}\mathcal{Y}_n$ according to the following rule:
$$
\lambda \cdot \mu = \sum{\nu} c^\nu_{\lambda, \mu} \nu, \lambda, \mu, \nu \in \mathcal{Y}n,
$$
where $c^\nu_{\lambda, \mu}$ are the Littlewood-Richardson coefficients,
then $\mathbb{Q}\mathcal{Y}_n$ is endowed with the structure of a commutative associative algebra over $\mathbb{Q}$, which will be isomorphic to the algebra $\Lambda_n$.
The  isomorphism $\varphi:\mathbb{Q}\mathcal{Y}_n  \longrightarrow \Lambda_n$ has the form  $\varphi(\lambda)=\boldsymbol{s}_\lambda$. Let us  introduce the following three operators on $\mathbb{Q}\mathcal{Y}_n$:
$$
\widehat D_\pm = \varphi^{-1} D_\pm \varphi, \widehat D_0 = \varphi^{-1} D_0 \varphi.
$$
It is evident that the operators $\widehat D_-, \widehat D_+, \widehat D_0$ defined in this way induce an $\mathfrak{sl}_2$ representation on $\mathbb{Q}\mathcal{Y}_n$, which is equivalent to the representation $\rho$ on $\Lambda_n$.

Consider the following three operators on $\mathcal{Y}_n$:

  \begin{align*}
 &\xi_-(\lambda)=\sum_{ \mu=\lambda-\square \in \mathcal{Y}_n}  \, \mu,\\
&\nabla_{\pm}(\lambda)=\sum_{\mu=\lambda \pm \square \in  \mathcal{Y}_n} c(\square) \, \mu,
\end{align*}

The following two propositions are reformulations of the results obtained earlier for symmetric polynomials in terms of Young diagrams.

\begin{te}\label{mt0} The triple of linear operators $(\widehat D_-,\widehat D_0,  \widehat D_+ )$ on $\mathbb{Q}\mathcal{Y}_n$, which acts on a partition $\lambda \in \mathbb{Q}\mathcal{Y}_n$   as follows:
\begin{align*}
 &\widehat  D_-(\lambda)=-( n \xi_-(\lambda)+\nabla_-(\lambda)),\\
&\widehat D_0(\lambda)=2 |\lambda| \lambda, \\
&\widehat D_+(\lambda)=\nabla_+(\lambda)
\end{align*}
defines an $\mathfrak{sl}_2$-action on $\mathbb{Q}\mathcal{Y}_n$.
\end{te}

\begin{te} The triple of linear operators $(\widetilde  D_-,\widetilde  D_0,  \widetilde  D_+ )$ on $\mathbb{Q}\mathcal{Y}^{(d)}_n$, which acts on a partition $\lambda \in  \mathcal{Y}^{(d)}_n$ as follows:
\begin{align*}
 &\widetilde{D}_-(\lambda) = -(n \xi_-(\lambda) + \nabla_-(\lambda)),\\
 &\widetilde  D_0(\lambda) = (2|\lambda|-n d)\lambda,\\
 &\widetilde D_+(\lambda) = \nabla_+(\lambda)-\ytableausetup{smalltableaux} \underbrace{\ydiagram{4}\cdots \ydiagram{1}}_{n \text{ boxes}} \cdot  \lambda,
\end{align*}
defines a finite-dimensional representation of $\mathfrak{sl}_2$  on $\mathbb{Q}\mathcal{Y}^{(d)}_n$.
\end{te}

We can also explicitly describe the decomposition of the $\mathfrak{sl}_2$-module $\mathbb{Q}\mathcal{Y}_n$ into irreducible submodules. First, we find the kernel of the operator $\widehat D_-$. To do this, we need to find the preimage of the elements $z_i$ under the isomorphism $\varphi$. We express the power sums $p_k$ in terms of Schur polynomials. From \cite[Theorem~7.17.1]{St}, we have that 

$$
p_k=s_{(k)}-s_{k-1,1}+s_{k-2,1,1}-s_{k-3,1,1,1}+\cdots+ (-1)^{k-1} s_{1,1,\ldots,1}, k=1,2,3,\ldots.
$$
Thus,
\begin{gather*}
\pi_k=\varphi^{-1}(p_k)=(k)-(k-1,1)+(k-2,1,1)-(k-3,1,1,1)+\cdots+ (-1)^{k-1} (\underbrace{1,1,\ldots,1}_{k \text{ boxes}})=\\
=\underbrace{\ydiagram{4}\cdots \ydiagram{1}}_{k \text{ boxes}}-\begin{ytableau}
\empty& \empty& \empty &  \none[\,...] &\empty\\
\empty& \none
\end{ytableau}+\begin{ytableau}
\empty& \empty& \empty &  \none[\,...] &\empty\\
\empty& \none\\
\empty& \none
\end{ytableau}+\cdots+(-1)^{k-1} \begin{ytableau}
\empty \\ \empty \\ \empty \\ \none \\\none[\vdots] \\ \phantom{\square}
\end{ytableau}\,.
\end{gather*}

As a vector space, $ \ker \widehat D_-$ is generated by the monomials $\zeta^\alpha=\zeta_2^{\alpha_1} \zeta_3^{\alpha_2} \cdots \zeta_n^{\alpha_{n-1}}$ where $\alpha \in \mathbb{N}^{n-1}$, where 
$$
\zeta_i=\varphi^{-1}(z_i)= \sum_{k=0}^{i-2} (-1)^k n^{i-k-1} {i \choose k} \pi_{i-k} \pi_1^k  +(i-1)(-1)^{i+1} \pi_1^i, i=2,3\ldots,n.
$$

Let $\widehat W(\zeta^\alpha)$ denote the vector space generated by the polynomials $\widehat D_+^k(\zeta^\alpha)$, $k \in \mathbb{N}$.

The following theorem is a reformulation of Theorem~\ref{dc} in terms of Young diagrams.

\begin{te}
The vector space $
\widehat W(\zeta^\alpha), 
$
is isomorphic to the standard infinite-dimensional lowest weight $\mathfrak{sl}_2$-module $W_{w(\alpha)}$ with the weight 
$$
	w(\alpha)=2 (2 \alpha_1+ 3 \alpha_2+\cdots+n \alpha_{n-1}).
	$$
For  $n>2$, we have a direct sum decomposition
$$
\mathbb{Q}\mathcal{Y}_n = \bigoplus_{\alpha \in \mathbb{N}^{n-1}}  \widehat W(\zeta^\alpha),
$$
and
$$
\mathbb{Q}\mathcal{Y}_n \cong \bigoplus_{i=0}^\infty c_i W_{2i}, 
$$
where $c_i$ is the number of natural solutions to the equation $2 \alpha_1+ 3 \alpha_2+\cdots+n \alpha_{n-1}=i$ and $W_{2i}$ is the standard infinite-dimensional $\mathfrak{sl}_2$-module of the type $([\circ])$ with weight $2i$.

\end{te}

From Theorem~\ref{fc}, it follows that the following decomposition holds:
$$
\mathbb{Q}\mathcal{Y}_n^{(d)} \cong \bigoplus\limits_{i=0}^{nd} c_n(d,i) V_i.
$$
Also, from Theorem~\ref{res}, it follows that the partitions $(1^i), i \leq d$ generate the standard finite-dimensional irreducible $\mathfrak{sl}_2$-module $V_d$.



\begin{thebibliography}{30}

\bibitem{HT}  Howe, R.,  Tan, E. C. (1992). Nonabelian harmonic analysis. Universitext. Springer-Verlag, New York, 10, 978-1.

\bibitem{G}
Gould, H. W. (1999). The Girard-Waring power sum formulas for symmetric functions, and Fibonacci sequences. Fibonacci Quarterly, 37, 135-140.

\bibitem{FH} Fulton, W.,  Harris, J. (2013). Representation theory: a first course (Vol. 129). Springer Science \& Business Media.
\bibitem{Essen}  Van den Essen, A. (2012). Polynomial Automorphisms: and the Jacobian Conjecture (Vol. 190). Birkh\"auser.
\bibitem{Essen1} van den Essen, A. (1995). Locally nilpotent derivations and their applications, III. Journal of Pure and Applied Algebra, 98(1), 15-23.
\bibitem{Now} Nowicki, A. (1994). Polynomial derivations and their rings of constants. Toru\'n: Uniwersytet Mikolaja Kopernika.

\bibitem{W23} Weigandt~A. Derivatives and Schubert Calculus. Talk at FPSAC conference 2023. URL: https://www.youtube.com/watch?v=fgzB7YjGOEc.

\bibitem{Gr} Grinberg, D., Korniichuk, N., Molokanov, K., \& Khomych, S. (2024). The diagonal derivative of a skew Schur polynomial. arXiv preprint arXiv:2402.14217.

\bibitem{Ernst}  Ernst, T. (2000). Generalized vandermonde determinants. Uppsala University. Department of Mathematics, Report 2000:6

\bibitem{KO}
Kamran, N.,  Olver, P. J. (1990). Lie algebras of differential operators and Lie-algebraic potentials. Journal of Mathematical Analysis and Applications, 145(2), 342-356.


\bibitem{Sp}
Springer, T. A. (1977). Invariant theory (Vol. 585). Springer.


\bibitem{St} Stanley, R.(2001) Enumerative Combinatorics. Volume 2. Cambridge University Press, Cambridge.
\bibitem{Pet} Petrov, L. (2013). \(\mathfrak{sl}(2)\) operators and Markov processes on branching graphs. Journal of Algebraic Combinatorics, 38(3), 663-720.

\bibitem{St88} Stanley, R. P. (1988). Differential posets. Journal of the American Mathematical Society, 1(4), 919-961.






\end{thebibliography}
\end{document}